\documentclass[english]{article}
\usepackage{lmodern}
\usepackage[T1]{fontenc}
\usepackage[utf8]{inputenc}

\usepackage{geometry}
\usepackage{extpfeil}
\usepackage[new]{old-arrows}
\usepackage{comment}
\usepackage[shortlabels]{enumitem}

\usepackage[dvipsnames,svgnames,table]{xcolor}

\usepackage[maxbibnames=99, style=numeric]{biblatex}
\AtBeginBibliography{\setlength{\emergencystretch}{2em}}

\usepackage{mathtools}
\usepackage{palatino}
\usepackage{mathpazo}
\usepackage{mathrsfs}

\usepackage{amsmath}
\usepackage{amsthm}
\usepackage{amssymb}
\usepackage{mathdots}
\usepackage{todonotes}
\usepackage{debate}
\usepackage{nicefrac}
\usepackage{graphicx}
\usepackage{caption}
\usepackage{float}
\DeclareFontFamily{U}{min}{}
\DeclareFontShape{U}{min}{m}{n}{<-> udmj30}{}
\newcommand\yo{\!\text{\usefont{U}{min}{m}{n}\symbol{'207}}\!}
\usepackage[all]{xy}
\usepackage[unicode=true,pdfusetitle,
 bookmarks=true,bookmarksnumbered=false,bookmarksopen=false,
 breaklinks=false,pdfborder={0 0 0},pdfborderstyle={},backref=false,colorlinks=true]{hyperref}

\usepackage{graphicx}

\hypersetup{
    colorlinks, linkcolor=OliveGreen,
    citecolor=OliveGreen, urlcolor=OliveGreen
}
 
\usepackage[nameinlink,capitalise,noabbrev]{cleveref}

\usepackage{bookmark}
\bookmarksetup{numbered}

\makeatletter

\RequirePackage{tikz-cd}
\RequirePackage{amssymb}
\usetikzlibrary{calc}
\usetikzlibrary{decorations.pathmorphing}

\tikzset{curve/.style={settings={#1},to path={(\tikztostart)
    .. controls ($(\tikztostart)!\pv{pos}!(\tikztotarget)!\pv{height}!270:(\tikztotarget)$)
    and ($(\tikztostart)!1-\pv{pos}!(\tikztotarget)!\pv{height}!270:(\tikztotarget)$)
    .. (\tikztotarget)\tikztonodes}},
    settings/.code={\tikzset{quiver/.cd,#1}
        \def\pv##1{\pgfkeysvalueof{/tikz/quiver/##1}}},
    quiver/.cd,pos/.initial=0.35,height/.initial=0}

\tikzset{tail reversed/.code={\pgfsetarrowsstart{tikzcd to}}}
\tikzset{2tail/.code={\pgfsetarrowsstart{Implies[reversed]}}}
\tikzset{2tail reversed/.code={\pgfsetarrowsstart{Implies}}}
\tikzset{no body/.style={/tikz/dash pattern=on 0 off 1mm}}

\newtheoremstyle{ctheorem}{}{}{}{}{\color{black}\bfseries}{}{ }{}
\theoremstyle{ctheorem}
\theoremstyle{definition}
\newtheorem{thm}{Theorem}[section]
\newtheorem{defn}[thm]{Definition}
\newtheorem{notation}[thm]{Notation}
\newtheorem{convention}[thm]{Convention}

\newtheorem{recollection}[thm]{Recollection}
\newtheorem{conjecture}[thm]{Conjecture}
\theoremstyle{definition}
\newtheorem{prop}[thm]{Proposition}
\theoremstyle{definition}
\newtheorem{lem}[thm]{Lemma}
\theoremstyle{definition}
\newtheorem{obsv}[thm]{Observation}
\theoremstyle{definition}
\newtheorem{rem}[thm]{Remark}
\theoremstyle{definition}
\newtheorem{example}[thm]{Example}
\theoremstyle{definition}
\newtheorem{cor}[thm]{Corollary}
\theoremstyle{definition}

\theoremstyle{definition}

\theoremstyle{definition}

\theoremstyle{definition}
\newtheorem{question}[thm]{Question}

\theoremstyle{definition}

\usepackage{mymacros}
\usepackage{citationscheme}
\author{Yuxuan Hu\thanks{Northwestern University.}}

\begin{document}
\title{Proper kernels in microlocal sheaf theory}
\maketitle
\begin{abstract}
Let $X$ and $Y$ be real analytic manifolds and let $\Lambda \subseteq T^*X$ and
$\Sigma \subseteq T^*Y$ be closed conic subanalytic singular isotropics.
Given a sheaf $K \in \Sh_{-\Lambda \times \Sigma}(X \times Y)$ microsupported
in $-\Lambda \times \Sigma$, consider the convolution functor
\(
(-) \ast K \colon \Sh_{\Lambda}(X) \rightarrow \Sh_{\Sigma}(Y)
\)
from sheaves microsupported in $\Lambda$ to sheaves microsupported in $\Sigma$.
We show that the convolution functor $(-) \ast K$ preserves compact objects
if and only if for each $x \in X$, the restriction
$K|_{\{x\} \times Y} \in \Sh_\Sigma(Y)$ is a compact object.
By a result of Kuo--Li \cite{kuo-li-duality},
the functor sending a sheaf kernel $K$ to the convolution functor $(-) \ast K$ is
an equivalence
between the category $\Sh_{-\Lambda \times \Sigma}(X \times Y)$ of sheaves
microsupported in $-\Lambda \times \Sigma$ and the category of cocontinuous functors
from $\Sh_\Lambda(X)$ to $\Sh_\Sigma(Y)$.
We therefore classify all cocontinuous functors that preserve compact objects
between the two categories.
Our approach is categorical and uses a small amount of geometry.
We introduce the notion of a proper object in a compactly
generated stable $\infty$-category and study its properties under
strongly cocontinuous localizations to obtain the result.
The main geometric input is the analysis of compact and proper objects of
the category of $P$-constructible sheaves for a triangulation $P$ of a manifold $Z$
via the exit-path category $\Exit(Z, P) \simeq P$.
Along the way, we show that a sheaf $F \in \Sh_\Lambda(X)$
is proper if and only if it has perfect stalks,
recovering a result of Ganatra--Pardon--Shende that extends earlier work of Nadler.
\end{abstract}
\tableofcontents
\section{Introduction}
\subsection{Motivation and background}
Let $X$ be a real analytic manifold and let $\Lambda \subseteq T^*X$ be a closed conic subanalytic singular isotropic.
The work of Kuo--Li \cite{kuo-li-duality} shows that the category
$\Sh_\Lambda(X)$ of sheaves microsupported in $\Lambda$ is a \emph{dualizable}
stable $\oo$-category (cf. \cref{defn:dualizable}).
More precisely, let $Y$ be a real analytic manifold and $\Sigma$ be a closed conic subanalytic singular isotropic.
It is shown there \cite[Theorems 1.1--1.3]{kuo-li-duality} that
the dual of $\Sh_\Lambda(X)$ is given by $\Sh_{-\Lambda}(X)$:
\[
	\Sh_\Lambda(X)^\vee \simeq \Sh_{-\Lambda}(X),
\]
the K\"unneth formula holds:
\[
	\Sh_\Lambda(X) \otimes \Sh_\Sigma(Y) \simeq \Sh_{\Lambda \times \Sigma}(X \times Y),
\]
and that the equivalence
\[
	\Sh_{-\Lambda\times \Sigma}(X\times Y) \simeq \Sh_\Lambda(X)^\vee \otimes \Sh_\Sigma(Y) \simeq \Fun^L(\Sh_\Lambda(X), \Sh_\Sigma(Y))
\]
is given by the assignment
\[
	K \mapsto (-) \ast K,
\]
where
\[
	(-) \ast K \coloneqq \pi_{2,!}(\pi_1^* (-) \otimes K)
\]
is the convolution functor. Here $\Fun^L(-,-)$ denotes the category of cocontinuous functors, and $\pi_1 \colon X\times Y \rightarrow X$ and
$\pi_2 \colon X \times Y \rightarrow Y$ are the projections.

In conclusion, cocontinuous functors between the sheaf categories
$\Sh_\Lambda(X)$ and $\Sh_\Sigma(Y)$ are classified by sheaf kernels
$K \in \Sh_{-\Lambda\times \Sigma}(X\times Y)$ on $X \times Y$ microsupported
in $-\Lambda \times \Sigma$.

On the other hand, previous works \cite{nadler-zaslow, nadler-microlocal-branes, gps3} have drawn parallels between
microlocal sheaf theory and the theory of
Fukaya categories in various flavors.
In particular, assuming $\Lambda \subseteq T^*X$ contains the zero section $0_X$,
it is shown in \cite[Theorem 1.1 and Remark 1.2]{gps3} that there is an equivalence
between the category of sheaves of $\ZZ$-modules microsupported in $\Lambda$ and
the Ind-completion of the perfect modules over the partially wrapped Fukaya category of $T^*X$
stopped at $-\Lambda_\infty$:
\[
	\Sh_\Lambda(X; \ZZ\mathhyphen\Mod) \simeq \Ind(\Perf \Wc(T^*X, -\Lambda_\infty)).
\]
Here $\Lambda_\infty$ is the projection of $\Lambda \setminus 0_X$ to the
cosphere bundle $S^*X$.
Taking compact objects on both sides,
we have
\[
	\Sh_\Lambda(X;\ZZ\mathhyphen\Mod)^\omega \simeq \Perf \Wc(T^*X,-\Lambda_\infty).
\]
Here $\Perf \Wc$ is the idempotent completion of the pretriangulated
envelope of the wrapped Fukaya category.
In \cite{nadler-wrapped}, Nadler first introduced the notion of
\emph{wrapped sheaves}, which are by definition compact objects in $\Sh_\Lambda(X)$.
In \cite{kuo-wrapped-sheaves}, Kuo showed that the category of wrapped sheaves
can also be realized geometrically, without alluding to compactness
in the ambient category.

Consequently, one can argue that studying compact objects
in $\Sh_\Lambda(X)$ is not a purely academic pursuit,
but is of intrinsic geometric interest from the viewpoint of Floer theory.
Given Liouville manifolds $M$ and $N$,
and an admissible Lagrangian correspondence $\Lc \subseteq M^{-} \times N$,
Gao \cite[Theorem 1.2]{gao-correspondence} showed that, if $\Lc \rightarrow N$ is proper,
under some genericity conditions,
there is an induced $A_\infty$-functor
\[
	\Theta_{\Lc} \colon \Wc(M) \rightarrow \Wc^{\mathrm{im}}(N),
\]
where $\Wc^{\mathrm{im}}(N)$ allows immersed Lagrangians with bounding cochains.
On objects the functor takes $L$ to the geometric composition
\[
	L \mapsto L \circ \Lc
\]
equipped with its canonical bounding cochain. A natural question is what the sheaf-theoretic incarnation
of the above functor is.
Note that the category
\[
	\Fun^\ex(\Sh_\Lambda(X)^\omega, \Sh_\Sigma(Y)^\omega)
\]
of exact functors is equivalent, by Ind-completion, to the full subcategory of
\[
	\Fun^L(\Sh_\Lambda(X), \Sh_\Sigma(Y)) \simeq \Sh_{-\Lambda \times \Sigma}(X \times Y)
\]
spanned by functors that preserve compact objects.
Therefore, an equivalent question to ask is:
\begin{question}
Under what conditions on the sheaf kernel
$K \in \Sh_{-\Lambda \times \Sigma}(X\times Y)$ does the convolution functor
\[
	(-) \ast K \colon \Sh_{\Lambda}(X) \rightarrow \Sh_{\Sigma}(Y)
\]
preserve compact objects?
\end{question}
Our main result classifies all such $K$ (see \cref{main-thm})
and implies the following assertion
from an early draft of Ganatra--Kuo--Li--Wu (see \cref{example:conj-ganatra}).

\begin{conjecture}[Ganatra--Kuo--Li--Wu]\label{conj-ganatra}
Let $Y$ be compact and let $K \in \Sh_{-\Lambda \times \Sigma}(X \times Y)$ be a
sheaf kernel. If $\ssupp(K) \rightarrow T^*Y$ is proper and $K$ has perfect stalks,
then $(-) \ast K$ preserves compact objects.
\end{conjecture}

\begin{rem}
Ganatra--Kuo--Li--Wu \cite{Ganatra-Kuo-Li-Wu} use geometric wrapping techniques
to study kernels obtained from constructible sheaves on $X \times Y$.
In their setting, the positive wrapping functor of \cite{kuo-wrapped-sheaves}
realizes the localization
$\iota^*_{-\Lambda \times \Sigma} \colon \Sh(X \times Y) \to \Sh_{-\Lambda \times \Sigma}(X \times Y)$.
\end{rem}
By our theorem, convolution with $K$ preserves compact objects precisely when
each fiber $K|_{\{x\}\times Y}$ is compact in $\Sh_\Sigma(Y)$.
When $K$ is obtained by wrapping a sheaf $L$, Ganatra--Kuo--Li--Wu give
geometric conditions on the support and microsupport of $L$ that guarantee
each fiber of $K$ is compact.
\begin{recollection}[{\cite[Theorem 5.20]{Ganatra-Kuo-Li-Wu}}]
Let $\Lambda\subseteq T^*X$ and $\Sigma\subseteq T^*Y$ be conic subanalytic Lagrangians,
and let $L$ be a real-constructible sheaf on $X\times Y$ with perfect stalks.
Assume that the projection $\supp(L)\to X$ is proper and that $\ssupp(L)$ meets
$-\Lambda\times\Sigma$ only in the zero section of $T^*(X\times Y)$. Set
\[
  \Gamma=\{(x,-\xi)\mid (x,\xi;y,\eta)\in\ssupp(L)
    \text{ for some }(y,\eta)\in\Sigma\}\subseteq T^*X.
\]
Assume finally that the link $\Gamma_\infty\subseteq S^*X$ is negatively
squeezable to $\Lambda_\infty$ with bounded energy \cite[Definition 5.10]{Ganatra-Kuo-Li-Wu}.
Such a squeezing is given by a negative contact isotopy that fixes $\Lambda_\infty$ and carries $\Gamma_\infty$
into arbitrarily small neighborhoods of a positive Reeb push-off of $\Lambda_\infty$,
with uniformly bounded Hamiltonian along the moving link.
Under these assumptions, convolution with the wrapped kernel $K=\iota^*_{-\Lambda\times\Sigma}(L)$
restricts to a functor
\[
  (-)\ast K\colon \Sh_\Lambda(X)^\omega\longrightarrow\Sh_\Sigma(Y)^\omega.
\]
\end{recollection}
In the special case $\ssupp^\infty(L)\subseteq T^*X\times(S^*Y\setminus\Sigma_\infty)$,
Ganatra--Kuo--Li--Wu \cite[\S~6.2]{Ganatra-Kuo-Li-Wu} identify each fiber
$K|_{\{x\}\times Y}$ with the wrapping of $L|_{\{x\}\times Y}$ in $\Sh_\Sigma(Y)$.
Each fiber is therefore compact, and \cref{main-thm} gives their conclusion.
\subsection{Main results and overview}
In this section, fix real analytic manifolds $X$ and $Y$
and closed conic subanalytic singular isotropics $\Lambda \subseteq T^*X$
and $\Sigma \subseteq T^*Y$.
Our main result is the following.
\begin{thm}[\cref{main-thm}]
Let $K \in \Sh_{-\Lambda\times \Sigma}(X\times Y)$ be a sheaf kernel.
The convolution functor
\[
	(-) \ast K\colon \Sh_\Lambda(X) \rightarrow \Sh_\Sigma(Y)
\]
preserves compact objects if and only if for every $x \in X$,
the restriction $K|_{\{x\} \times Y} \in \Sh_{\Sigma}(Y)$ is a compact object.
Consequently, there is an equivalence
\[
	\Pc \simeq \Fun^{\ex}(\Sh_\Lambda(X)^\omega,\Sh_\Sigma(Y)^\omega),
\]
where $\Pc \subseteq \Sh_{-\Lambda \times \Sigma}(X \times Y)$
is the full subcategory spanned by such sheaf kernels.
\end{thm}
From this, we can deduce a sufficient condition.
\begin{cor}[\cref{perfect-stalk-and-compact-support}]
Let $K \in \Sh_{-\Lambda \times \Sigma}(X \times Y)$ have perfect stalks.
If its restriction to $\{x\} \times Y$ has compact support for every $x \in X$,
then convolution with $K$ preserves compact objects.
\end{cor}
\begin{rem}\label{example:conj-ganatra}
If $Y$ is compact, then the closed support $\supp(K|_{\{x\} \times Y})$ is compact.
In this case, if $K$ has perfect stalks, then convolution with $K$ preserves
compact objects. In particular, \cref{conj-ganatra} is true.
\end{rem}
The core of our argument is the notion of a
\emph{proper} object
in a compactly generated stable $\oo$-category.
\begin{defn}[\cref{defn:proper-object}]
Let $\Cc$ be a presentable stable $\infty$-category.
We say $c\in\Cc$ is \emph{proper} if the functor
\[
\map_\Cc(-,c) \colon \Cc^{\omega,\op} \to \Sp
\]
factors through $\Sp^\omega \subseteq \Sp$.
\end{defn}
The following result on proper objects in $\Sh_\Lambda(X)$ can also be
seen as a special case of the main theorem, using the duality above and taking
$Y$ to be a point with $\Sigma = T^*Y$.
\begin{cor}[\cref{proper-objects-in-shv}]
A sheaf $F \in \Sh_\Lambda(X)$ is proper if and only if $F$ has perfect stalks.
\end{cor}
Nadler proved a microlocal version using arborealization
\cite[Theorem 3.21]{nadler-wrapped}. Ganatra--Pardon--Shende gave a direct proof
for sheaves of modules in \cite[Corollary 4.24]{gps3}.
The main point of this paper is that,
once we distill the essential ideas in its proof
in categorical terms, the argument can be
further simplified and adapted to a relative setting,
allowing us to prove the main theorem of this paper.
\subsection{Notations and conventions}
For the sake of brevity and clarity,
we will work exclusively with sheaves of spectra, unless otherwise specified.
\begin{rem}
All our arguments work mutatis mutandis if
we replace the $\infty$-category $\Sp$ of spectra with
a compactly generated rigid symmetric monoidal stable $\infty$-category $\Vc$,
using $\Vc$-linear functors and tensor products over $\Vc$.
Much of the theory of $\Vc$-enriched categories is developed in
\cite{gepner-haugseng, hinich18, hinich21, heine-equivalence, benmoshe-naturality}.
For a quick review on the theory of presentable and dualizable categories
in the enriched setting directly applicable to this paper, see \cite[\S 1]{ramzi2024dualizable}.
\end{rem}
\begin{notation}
Let $\Cc$ be a stable $\oo$-category.
Throughout this paper,
$\map_\Cc(-,-)$ denotes the mapping spectrum,
while $\Map_\Cc(-,-)$ denotes the mapping space.
\end{notation}
\subsection*{Acknowledgements}
I am grateful to Haosen Wu for bringing this problem to my attention and sharing an early draft of the joint work \cite{Ganatra-Kuo-Li-Wu}, and to Wenyuan Li for answering many questions about that work.
This project began as an attempt to give a categorical proof of a geometric result pursued in \cite{Ganatra-Kuo-Li-Wu}.
While the theorem we prove here is more general than originally anticipated, it addresses a different (though related) question than the one studied there.
I thank Qingyuan Bai, Peter Haine, and Gus Schrader for reading an early draft of this paper, and Peter Haine in particular for very detailed comments.
I was supported by the NSF grant DMS 2302624.

\subsection*{AI Declaration}
An earlier version of this paper (\href{https://arxiv.org/abs/2511.02677v1}{arXiv:2511.02677v1}) was written without AI assistance.
For this revision, OpenAI Codex (GPT-6.2 Astra and Sol) was used to correct typos, reorganizing the presentation of existing arguments, and update the comparison with Ganatra-Kuo-Li-Wu.

\section{Preliminaries}
\subsection{Exodromy}
Consider a topological space $X$ with a stratification $P$,
or simply a \emph{stratified space} $(X,P)$.
\begin{recollection}
A \emph{stratified space} $(X,P)$ is a topological space $X$,
together with a continuous map $X \rightarrow P$, where $P$ is a poset
with the Alexandroff topology whose open sets are the upward-closed subsets.
\end{recollection}
The \emph{exodromy equivalence}
(see \cite{treumann2009exit, HA, lejay2021constructible, porta-teyssier, haine2024exodromy})
states that under suitable assumptions
on $(X,P)$ and $\Vc$, the category
\[
	\Cons_P(X;\Vc)
\]
of $P$-constructible sheaves\footnote{Technically, to state the most general result,
one has to consider hyper-constructible hypersheaves.
However, the distinction between hyper-constructible hypersheaves
and constructible sheaves disappears when everything is hypercomplete,
as is the case if $X$ is a finite-dimensional manifold and all the strata
are submanifolds.}
on $X$ valued in $\Vc$
is equivalent to the category
\[
	\Fun(\Exit(X,P), \Vc)
\]
of functors from the \emph{exit-path $\oo$-category}
$\Exit(X,P)$ to $\Vc$.
While the construction of $\Exit(X,P)$ is quite involved in general,
the only result we need in this paper is the following.
\begin{prop}
Let $P$ be the face poset of a locally finite triangulation of a manifold $X$. Then $\Exit(X,P) \simeq P$.
\end{prop}
\begin{proof}
This is a special case of \HA{Theorem}{A.6.10}.
\end{proof}
\begin{cor}\label{cons-p-fun-p}
Let $P$ be the face poset of a locally finite triangulation of a manifold $X$.
Then $\Cons_P(X) \simeq \Fun(P, \Sp)$.
\end{cor}

\subsection{Dualizable stable $\infty$-categories}
We give a quick recap of the theory of dualizable stable $\oo$-categories.

The $\infty$-category $\Prlst$ of presentable stable $\infty$-categories and left adjoints admits
a closed symmetric monoidal structure given by tensor product in $\Prl$.
The internal hom is given by $[\Cc, \Dc] = \Fun^L(\Cc, \Dc)$.

The category $\Catperf$ consists of small idempotent-complete stable
$\infty$-categories and exact functors. In its symmetric monoidal structure,
the tensor product $\Cc_0 \otimes \Dc_0$
classifies bi-exact functors out of $\Cc_0 \times \Dc_0$.
It is again a closed monoidal category with internal homs given by
$[\Cc_0, \Dc_0] = \Fun^{\ex}(\Cc_0, \Dc_0)$.

The ind-construction
\[
\Ind \colon \Catperf \to \Prlst
\]
identifies $\Catperf$ with a subcategory of $\Prlst$, with
its essential image spanned by compactly generated stable $\infty$-categories
and left adjoints that preserve compact objects.

Moreover, $\Ind$ is a symmetric monoidal functor.
The ind construction and taking compact objects give inverse equivalences
\[
	\Ind \colon \Catperf \leftrightarrows {\Pr}^{L,\omega}_\st \colon (-)^\omega,
\]
where $\Pr^{L,\omega}_\st$ is the $\oo$-category of compactly generated
stable $\oo$-categories and left adjoints that preserve compact objects.
\begin{defn}\label{defn:dualizable}
A presentable stable $\infty$-category $\Cc$ is \emph{dualizable} if
there exists $\Cc^\vee$ so that the functor
\[
- \otimes \Cc \colon \Prlst \to \Prlst
\]
is left adjoint to
\[
- \otimes \Cc^\vee \colon \Prlst \to \Prlst.
\]
\end{defn}
\begin{example}[{{\SAG{Proposition}{D.7.2.3}}}]
If $\Cc$ is a compactly generated stable $\oo$-category, then $\Cc$ is
dualizable. The dual is given by $\Cc^\vee \simeq \Ind(\Cc^{\omega,\op})$,
and the evaluation map $\ev: \Cc^\vee \otimes \Cc \to \Sp$ is given by
the left Kan extension of
\[
  \map_\Cc(-, -) \colon \Cc^{\omega,\op} \otimes \Cc^\omega \to \Sp
\]
along $\Cc^{\omega,\op} \otimes \Cc^\omega \to \Cc^\vee \otimes \Cc \simeq \Ind(\Cc^{\omega,\op} \otimes \Cc^\omega)$.
\end{example}
\begin{recollection}
A left adjoint functor between presentable $\oo$-categories is called
an \emph{internal left adjoint} in $\Pr^L$, if it is the left adjoint of
an adjunction in the $(\infty,2)$-category $\mathbf{Pr}^L$ of presentable $\oo$-categories and left adjoints.
Being an internal left adjoint is equivalent to being a \emph{strongly cocontinuous functor}: a functor whose right adjoint admits a further right adjoint.
\end{recollection}
\begin{rem}\label{internal-left-adjoint-between-compactly-generated}
A left adjoint functor between compactly generated presentable stable $\infty$-categories
is an internal left adjoint if and only if it preserves compact objects.
Therefore, $\Catperf \simeq \Pr^{L,\omega}_\st$ is equivalent to
the $\oo$-category of compactly generated stable $\oo$-categories and internal
left adjoints.
\end{rem}
\begin{defn}
The $\infty$-category $\Pr^{\dual}_\st$ of dualizable stable $\infty$-categories
is the subcategory of $\Prlst$ with dualizable objects
and internal left adjoints as morphisms.
\end{defn}
\begin{rem}
By \cref{internal-left-adjoint-between-compactly-generated}, $\Catperf \simeq \Pr^{L,\omega}_\st$ identifies with the full subcategory of $\Pr^\dual_\st$ spanned
by compactly generated stable $\oo$-categories.
\end{rem}
\begin{thm}[{{\SAG{Proposition}{D.7.3.1}}}]
  A presentable stable $\infty$-category is dualizable
  if and only if it is a retract of a compactly
  generated one in $\Prlst$.
\end{thm}
\section{Compactness and properness}
In this section, we discuss \emph{compact} and \emph{proper} objects
in presentable stable $\infty$-categories.
We will focus on how these objects behave under
the inclusion of reflective and coreflective subcategories.
\begin{recollection}
Let $\mathcal{C}$ be a presentable $\infty$-category.
An object $c \in \mathcal{C}$ is \emph{compact} if the functor
$\Map_{\Cc}(c, -) \colon \mathcal{C} \to \Spc$ commutes with filtered colimits.
If $\Cc$ is stable, this is equivalent to requiring that
$\map_{\Cc}(c, -) \colon \Cc \to \Sp$ commutes with filtered colimits.
\end{recollection}
\begin{defn}\label{defn:proper-object}
Let $\Cc$ be a presentable stable $\infty$-category.
We say $c\in\Cc$ is \emph{proper} if the functor
\[
\map_\Cc(-,c) \colon \Cc^{\omega,\op} \to \Sp
\]
factors through $\Sp^\omega \subseteq \Sp$.
Denote by $\Cc^{\text{pr}} \subseteq \Cc$ the full subcategory of proper objects.
\end{defn}
\begin{rem}
If $\Cc$ is compactly generated and every compact object in $\Cc$ is proper,
then $\Cc$ is a proper stable $\infty$-category in the sense of \SAG{Definition}{11.1.0.1}.
\end{rem}
\begin{prop}
Let $\Cc$ be a presentable stable $\oo$-category.
Then $\Cc^{\text{pr}}$ is an idempotent-complete stable subcategory of $\Cc$.
\end{prop}
\begin{proof}
Note that $\Sp^\omega \subseteq \Sp$ is closed under finite (co)limits and retracts,
and for any $x \in \Cc^\omega$,
the functor $\map_\Cc(x,-) \colon \Cc \rightarrow \Sp$ preserves finite (co)limits and retracts.
It follows that $\Cc^\pr$ is closed under finite (co)limits and retracts.
\end{proof}
\begin{prop}
Let $\Cc$ be a compactly generated stable $\oo$-category.
The spectral Yoneda embedding
\[
	c  \mapsto \map_\Cc(-,c)
\]
restricts to an equivalence
\[
	\Cc^\pr \xrightarrow{\simeq} \Fun^\ex(\Cc^{\omega,\op}, \Sp^\omega).
\]
Here $\Fun^\ex(-,-)$ denotes the category of exact
functors between two stable $\oo$-categories.
\end{prop}
\begin{proof}
Since $\Cc$ is compactly generated, there are equivalences
\[
	\Cc \simeq \Ind(\Cc^\omega) \simeq \Fun^{\lex}(\Cc^{\omega,\op},\Spc)
	\simeq \Fun^{\ex}(\Cc^{\omega,\op},\Sp).
\]
Here $\Fun^{\lex}(-,-)$ denotes the category of left exact\footnote{Recall that a functor is called \emph{left exact} if it preserves finite limits.} functors.
By definition, the composite functor is the spectral Yoneda embedding.
Restricting to proper objects gives the desired equivalence.
\end{proof}
\begin{cor}
If $\Cc$ is a compactly generated stable $\oo$-category,
then $\Cc^\pr$ is a small $\oo$-category.
\qed
\end{cor}
Compact and proper objects behave in a very controllable way under
\emph{strongly cocontinuous localizations}, which we now introduce.
\begin{recollection}
A functor between presentable $\infty$-categories is called \emph{strongly cocontinuous},
if its right adjoint admits a further right adjoint.
Alternatively, it is an internal left adjoint in $\Pr^L$.
\end{recollection}
\begin{obsv}
Let $\Cc$ be a presentable $\oo$-category, and $\iota_* \colon \Dc \hookrightarrow \Cc$
be a full subcategory.
Suppose that $\Dc$ is closed under small limits and colimits.
By the $\oo$-categorical reflection theorem \cite{infinity-categorical-reflection},
the inclusion $\iota_*$ participates in a triple adjunction:
\[\begin{tikzcd}[cramped]
	\Dc & \Cc
	\arrow["{\iota_*}"{description}, hook, from=1-1, to=1-2]
	\arrow["{\iota^\flat}", shift left=3, from=1-2, to=1-1]
	\arrow["{\iota^*}"', shift right=3, from=1-2, to=1-1]
\end{tikzcd}.\]
In particular, the localization functor $\iota^*$ is strongly cocontinuous.
\end{obsv}
\begin{notation}\label{bireflective-subcategory}
In the above situation, we say $\iota^* \colon \Cc \rightarrow \Dc$
is a \emph{strongly cocontinuous localization}, and $\iota_* \colon \Dc \hookrightarrow \Cc$
is a \emph{bi-reflective subcategory}. In this case, we always use
$\iota^* \dashv \iota_* \dashv \iota^\flat$ to refer to the adjoint triple.
\end{notation}
\begin{prop}
\label{localization-compactness}
Let $\Cc$ be a compactly generated presentable $\oo$-category.
Suppose $\iota_* \colon \Dc \hookrightarrow \Cc$ is a bi-reflective subcategory.
Then:
\begin{enumerate}[(1)]
   \item The inclusion $\iota_*$ detects compact objects: if $\iota_* d$ is compact in $\Cc$, then $d$ is compact in $\Dc$.
   \item The $\infty$-category $\Dc$ is compactly generated by $\iota^* \Cc^\omega$.
   \item The category $\Dc^\omega$ is the smallest replete subcategory of $\Dc$
	   containing $\iota^* \Cc^\omega$ and closed under finite colimits
	   and retracts.
\end{enumerate}
\end{prop}
\begin{proof}
To prove point (1), note that there is a natural equivalence
\[
\Map_\Dc(d, -) \simeq \Map_\Cc(\iota_*d, \iota_* -)
\]
and that $\iota_*$ preserves all colimits, in particular filtered colimits.
Therefore if $\iota_*d$ is compact in $\Cc$, then $\Map_\Dc(d, -)$ commutes with filtered colimits, and hence $d$ is compact in $\Dc$.

To prove points (2) and (3), first note that $\iota^*$ preserves compact objects:
indeed, its right adjoint $\iota_*$ preserves all colimits, in particular filtered colimits \HTT{Proposition}{5.5.7.2}.

Now let $\Dc'$ be the full subcategory of $\Dc$ generated under colimits
by $\iota^*\Cc^\omega$.
Consider the full subcategory
\[
  (\iota^*)^{-1} \Dc' \defeq \{ c \in \Cc \mid \iota^* c \in \Dc' \} \subseteq \Cc.
\]
Since $\iota^*$ preserves colimits, it is immediate that $(\iota^*)^{-1}(\Dc')$ is closed under colimits.
However, by definition of $\Dc'$ we have
\[
\Cc^\omega \subseteq (\iota^*)^{-1}\iota^* \Cc^\omega \subseteq (\iota^*)^{-1}\Dc'.
\]
As $\Cc^\omega$ generates $\Cc$ under colimits, we must have $(\iota^*)^{-1}\Dc' = \Cc$.
Therefore, we obtain $\Dc' = \Dc$: for any $d \in \Dc$, we have
$d \simeq \iota^* \iota_* d \in \iota^* (\Cc) = \iota^*(\iota^*)^{-1}\Dc' \subseteq \Dc'$.
Consequently, $\iota^*\Cc^\omega$ generates $\Dc$ under colimits.
Point (3) follows from (2) because a collection of compact generators under
colimits generates the compact objects under finite colimits and retracts
(cf. \HTT{Proposition}{5.3.5.11} and \HTT{Lemma}{5.4.2.4}).
\end{proof}
We use this description of $\Dc^\omega$ to test properness after localization
in \cref{v-properness-localization}, where we allow coefficients in a compactly
generated stable $\oo$-category.

\section{\texorpdfstring{$\Vc$}{V}-properness}
Throughout this section, let $\Vc$ be a compactly generated stable $\oo$-category.
\begin{notation}
Let $\Cc$ be a dualizable stable $\oo$-category.
Denote by
\[
	e \colon \Cc \otimes \Vc \xrightarrow{\simeq} \Fun^L(\Cc^\vee, \Vc)
\]
the evaluation functor.
\end{notation}
\begin{defn}
Let $\Cc$ be a compactly generated stable $\oo$-category.
An object $F \in \Cc \otimes \Vc$ is \emph{$\Vc$-proper} if the functor
\[
	e(F) \colon \Cc^\vee \rightarrow \Vc
\]
sends $(\Cc^\vee)^\omega \simeq \Cc^{\omega,\op}$ to $\Vc^\omega$.
\end{defn}
\begin{rem}
An object $x \in \Cc$ is proper in the sense of \cref{defn:proper-object}
precisely if $x \in \Cc \simeq \Cc \otimes \Sp$ is $\Sp$-proper.
\end{rem}
\begin{example}
Let $x \in \Cc$ be a proper object and $v \in \Vc$ a compact object.
Then $x \boxtimes v \in \Cc \otimes \Vc$ is $\Vc$-proper.
Here $- \boxtimes - \colon \Cc \times \Vc \rightarrow \Cc \otimes \Vc$ is
the universal bi-cocontinuous functor.
\end{example}
\begin{obsv}
Let $\Cc$ be a compactly generated stable $\oo$-category.
Suppose $\iota_* \colon \Dc \hookrightarrow \Cc$ is a bi-reflective subcategory (cf. \cref{bireflective-subcategory}):
  \[\begin{tikzcd}[cramped]
    \Dc & \Cc
    \arrow["{\iota_*}"{description}, hook, from=1-1, to=1-2]
    \arrow["{\iota^\flat}", shift left=3, from=1-2, to=1-1]
    \arrow["{\iota^*}"', shift right=3, from=1-2, to=1-1]
  \end{tikzcd}.\]
Since $\iota^* \dashv \iota_*$ is an adjunction internal to $\Pr^L$,
tensoring with $\Vc$ gives rise to a bi-reflective subcategory
\[\begin{tikzcd}[cramped]
	{\Dc \otimes \Vc} & {\Cc \otimes \Vc}
	\arrow["{\iota_*^\Vc}"{description}, hook, from=1-1, to=1-2]
	\arrow["{\iota^\flat_\Vc}", shift left=3, from=1-2, to=1-1]
	\arrow["{\iota^*_\Vc}"', shift right=3, from=1-2, to=1-1].
\end{tikzcd}\]
Here $\iota^*_\Vc = \iota^* \otimes \Vc$, $\iota_*^\Vc = \iota_* \otimes \Vc$, and $\iota^\flat_\Vc$ is the right adjoint to $\iota_* \otimes \Vc$.
\end{obsv}
\begin{prop}
\label{v-properness-localization}
In the above situation, an object $F \in \Dc \otimes \Vc$ is $\Vc$-proper
if and only if $\iota_*^\Vc(F)$ is $\Vc$-proper.
\end{prop}
\begin{proof}
By \HA{Proposition}{4.8.1.17},
there is a commutative diagram
\[\begin{tikzcd}[cramped]
	{} & \Prr && \Prr \\
	& {\Pr^{L,\op}} && {\Pr^R}
	\arrow["{- \otimes \Vc}", from=1-2, to=1-4]
	\arrow["\simeq"', from=1-2, to=2-2]
	\arrow[equals, from=1-4, to=2-4]
	\arrow["{\Fun^{R}(-^\op,\Vc)}", from=2-2, to=2-4]
\end{tikzcd}.\]
For compactly generated stable $\oo$-categories and right adjoints of
compactness-preserving left adjoints, we have
\[
	- \otimes \Vc \simeq \Fun^R(-^\op,\Vc) \simeq \Fun^\ex((-)^{\omega,\op},\Vc).
\]
Under this equivalence, the right adjoint $\iota_*^\Vc$ corresponds to
\[
	- \circ (\iota^*|_{\Cc^\omega})^\op \colon \Fun^{\ex}(\Dc^{\omega,\op},\Vc) \rightarrow \Fun^{\ex}(\Cc^{\omega,\op},\Vc).
\]
Suppose $F \colon \Dc^{\omega,\op} \rightarrow \Vc$ factors through $\Vc^\omega$.
It follows that $\iota_*^\Vc(F)$ also factors through $\Vc^\omega$.
On the other hand, suppose $\iota_*^\Vc(F)$ factors through $\Vc^\omega$,
and let $\mathcal{K} \subseteq \Dc^\omega$ be the full subcategory spanned
by objects $d \in \Dc^\omega$ such that $F(d) \in \Vc^\omega$.
By assumption on $\iota_*^\Vc(F)$, we have
$F(\iota^*c) = \iota_*^\Vc(F)(c) \in \Vc^\omega$,
and therefore $\iota^*\Cc^\omega \subseteq \mathcal{K}$.
Because $\mathcal{K}$ is closed under finite colimits and retracts,
\cref{localization-compactness} shows that $\mathcal{K} = \Dc^\omega$.
\end{proof}
\begin{cor}\label{properness-localization}
Let $\Cc$ be a compactly generated stable $\oo$-category.
Let $\iota_* \colon \Dc \hookrightarrow \Cc$ be a bi-reflective subcategory.
Then $d \in \Dc$ is proper in $\Dc$
if and only if $\iota_*d \in \Cc$ is proper in $\Cc$.
\end{cor}
\begin{proof}
Take $\Vc=\Sp$ in \cref{v-properness-localization}. Under the identification
$\Dc\otimes\Sp\simeq\Dc$, $\Sp$-properness is ordinary properness and
$\iota_*^\Sp$ identifies with $\iota_*$.
\end{proof}
\section{Proper objects in functor categories}
Let $P$ be the face poset of a locally finite triangulation of a manifold $X$.
By the exodromy equivalence (cf. \cref{cons-p-fun-p}), we have
\[
\Cons_P(X) \simeq \Fun(P, \Sp).
\]
In this section, we study proper objects in $\Cons_P(X) \simeq \Fun(P,\Sp)$.
First we note that, if $P$ is a triangulation of $X$,
then $P$ as a poset is \emph{locally finite}.
\begin{defn}
A poset $P$ is \emph{locally finite} if for every $p \in P$, the poset $P_{p/}$ is finite.
\end{defn}
\begin{prop}[{\cite[Lemma 4.4.10]{bh25}}]\label{compact-objects-in-fun-p}
Let $P$ be a locally finite poset.
Then $F$ is compact in $\Fun(P, \Sp)$ if and only if
$F$ is finitely supported and $F(p)$ is a finite spectrum for every $p \in P$.
\end{prop}
\begin{thm}\label{v-proper-objects-in-fun-p}
Let $P$ be a locally finite poset. Then $F \in \Fun(P,\Sp) \otimes \Vc \simeq \Fun(P,\Vc)$ is $\Vc$-proper if and only if $F(p) \in \Vc^\omega$ for every $p \in P$.
\end{thm}
\begin{proof}
Write
\begin{align*}
	\yo \colon P^\op & \rightarrow \Fun(P, \Spc) \rightarrow \Fun(P, \Sp) \\
	p & \mapsto \Map_P(p, -) \mapsto \Sigma^\infty_+\Map_P(p, -)
\end{align*}
for the stable (co)Yoneda embedding.

By \cite[Proposition 2.2.3]{LurieRotation}, $\Fun(P,\Sp)$ is compactly generated
by the collection $\{\yo(p)\}_{p \in P}$.
Recall that the equivalence
\[
	\Fun(P, \Sp) \otimes \Vc \simeq \Fun(P, \Vc)
\]
can be obtained as the composite
\[
	\Fun(P, \Sp) \otimes \Vc \simeq \Fun^R(\Fun(P,\Sp)^\op, \Vc)
	\simeq \Fun^R(\Fun(P,\Spc)^\op, \Vc) \simeq \Fun(P, \Vc).
\]
It follows that under the above equivalence,
the continuous functor
\[
	\tilde{e}(F) \colon \Fun(P,\Sp)^\op \rightarrow \Vc
\]
classified by an object $F \in \Fun(P,\Vc)$
sends $\yo(p)$ to $F(p) \in \Vc$.

On the other hand, unwinding definitions
(cf. the proof of \SAG{Proposition}{D.7.2.3}),
the evaluation
\[
	e(F) \colon \Fun(P,\Sp)^\vee \rightarrow \Vc
\]
is the ind extension of
\[
	\tilde{e}(F)|_{\Fun(P,\Sp)^{\omega,\op}} \colon \Fun(P,\Sp)^{\omega,\op} \rightarrow \Vc.
\]
Since $\yo(p)$ is compact in $\Fun(P,\Sp)$,
it follows that
\[
	e(F)(\yo(p)) \simeq \tilde{e}(F)(\yo(p)) \simeq F(p).
\]
The objects $\yo(p)$ generate $\Fun(P,\Sp)^\omega$ under finite limits,
finite colimits, and retracts. Since $e(F)$ is exact and these operations
preserve compact objects in $\Vc$, the result follows.
\end{proof}
\begin{cor}\label{properness-of-fun-p}
  If $P$ is a locally finite poset, then $\Fun(P, \Sp)$ is a proper stable $\infty$-category:
  every compact object is proper.
  The compact objects are precisely those functors that factor through $\Sp^\omega$
  and are supported on a finite subset of $P$.
  \qed
\end{cor}

\section{Proper objects in sheaf categories}
We prove the main theorem by reducing to the $P$-constructible sheaves
studied in the previous section.

To this end, we first recall some facts on the geometry of stratifications
and microlocal sheaf theory.
\begin{notation}
Let $P$ be a $C^1$-stratification of a $C^1$-manifold $X$.
We write
\[
	N^*P \coloneqq \cup_{p \in P} N^*X_p \subseteq T^*X
\]
for the union of the conormals of the strata in $X$.
\end{notation}
\begin{recollection}[{\cite[Lemma 2.12]{kuo-li-duality}, \cite{czapla2012definable}}]
Let $X$ be a real analytic manifold and $\Lambda \subseteq T^*X$ be a closed conic subanalytic singular isotropic.
There is a locally finite $C^1$ Whitney triangulation $P$ of $X$ such that $\Lambda \subseteq N^*P$.
\end{recollection}
\begin{recollection}[{{\cite[Proposition 8.4.1]{ks90}, \cite[Proposition 2.7]{kuo-li-duality}}}]
Let $P$ be a $C^1$ Whitney stratification of $X$.
Then
\[
	\Cons_P(X) \simeq \Sh_{N^*P}(X).
\]
\end{recollection}
\begin{recollection}[{\cite[\S 2.3]{kuo-li-duality}}]
Let $\Lambda \subseteq T^*X$ be a closed conic subanalytic singular isotropic.
The $\oo$-category $\Sh_\Lambda(X)$ is closed under limits and colimits in $\Sh(X)$.
\end{recollection}
The results mentioned so far can be combined and summarized as follows.
\begin{obsv}
Let $X$ be a real analytic manifold and $\Lambda \subseteq T^*X$ a closed conic
subanalytic singular isotropic.
Then there exists a locally finite $C^1$ Whitney triangulation $P$ of $X$, so that $\Lambda \subseteq N^*P$.
And the inclusion
\[
	\iota_{\Lambda, P, *} \colon \Sh_\Lambda(X) \hookrightarrow \Sh_{N^*P}(X) \simeq \Cons_P(X) \simeq \Fun(P, \Sp)
\]
exhibits $\Sh_\Lambda(X)$ as a bi-reflective subcategory of $\Cons_P(X)$
(cf. \cref{bireflective-subcategory}).
In particular, it participates
in an adjoint triple
$\iota^*_{\Lambda,P} \dashv \iota_{\Lambda,P,*} \dashv \iota_{\Lambda,P}^\flat$.
\end{obsv}
For the rest of this section, fix $X$, $\Lambda$, and $P$ as above.
By \cref{localization-compactness}, $\Sh_\Lambda(X)$ is compactly generated.
\begin{cor}\label{compacts-in-shv}
If $F \in \Sh_\Lambda(X)$ is compactly supported and has perfect stalks,
then $F$ is compact.

\end{cor}
\begin{proof}
The compact support of $F$ meets only finitely many strata of the locally
finite triangulation $P$. Thus by \cref{compact-objects-in-fun-p}, $F$ is compact in $\Cons_P(X)\simeq \Fun(P,\Sp)$ and
thus also compact in $\Sh_\Lambda(X)$ by \cref{localization-compactness}.
\end{proof}
To test properness with coefficients, we first identify the stalk functors
under the exodromy equivalence.
\begin{lem}\label{stalk-in-cons-p}
	Let $(X, P)$ be an exodromic stratified space, and $\Vc$ a dualizable stable $\oo$-category\footnote{Here the dualizability condition is assumed to ensure that we have\[
	\Cons_P(X;\Vc) \simeq \Cons_P(X;\Spc) \otimes \Vc \simeq \Fun(\Exit(X,P),\Vc).\] See {\cite[\S 4]{haine2024exodromy}} for detailed discussions on the exodromy equivalence with coefficients.}
Then the stalk functor at $x$
\[
(-)_x \colon \Cons_P(X;\Vc) \rightarrow \Vc
\]
is canonically equivalent to
the evaluation functor
\[
\ev_x \colon \Fun(\Exit(X, P), \Vc) \rightarrow \Vc
\]
at $x \colon [0] \rightarrow \Exit(X, P)$.
\end{lem}
\begin{proof}
Taking stalks at $x$ is by definition the pullback along
$x: \ast \rightarrow X$, \emph{i.e.}
\[
F_x \simeq x^* F \in \Sh(\ast; \Vc) \simeq \Vc.
\]
By the functoriality of the exodromy equivalence,
this is equivalent to evaluation at $x \colon [0] \rightarrow \Exit(X, P)$.
\footnote{Here we distinguish the topological space $\ast$ consisting of a single point
from its homotopy type $[0]$.}
\end{proof}

\begin{thm}\label{v-proper-object-in-shv}
Let $\Vc$ be a compactly generated stable $\oo$-category.
An object $F \in \Sh_\Lambda(X) \otimes \Vc$ is $\Vc$-proper if and only
if viewed as an object in $\Sh(X;\Vc)$, it has stalks valued in $\Vc^\omega$.
\end{thm}
\begin{proof}
By \cref{v-properness-localization}, $F \in \Sh_\Lambda(X) \otimes \Vc$ is
$\Vc$-proper if and only if
\[
	\iota_{\Lambda, P, *}^\Vc (F) \in \Cons_P(X) \otimes \Vc \simeq \Fun(P,\Sp) \otimes \Vc \simeq \Fun(P, \Vc)
\]
is $\Vc$-proper.
By \cref{v-proper-objects-in-fun-p},
this is equivalent to $F(p) \in \Vc^\omega$ for every $p \in P$,
which in turn is equivalent to $F$ having stalks valued
in $\Vc^\omega$ when viewed as an object in $\Sh(X;\Vc)$
by \cref{stalk-in-cons-p}.
\end{proof}
\begin{cor}\label{proper-objects-in-shv}
A sheaf $F \in \Sh_\Lambda(X)$ is proper if and only if $F$ has perfect stalks.
\end{cor}
\begin{proof}
Apply \cref{v-proper-object-in-shv} with $\Vc=\Sp$. The compact objects of
$\Sp$ are the perfect spectra.
\end{proof}
Combining \cref{v-proper-object-in-shv} with the following results of
Kuo--Li \cite{kuo-li-duality} gives our main theorem.
\begin{recollection}[{{\cite[Theorems 1.1--1.3]{kuo-li-duality}}}]
Let $X$ and $Y$ be real analytic manifolds, $\Lambda \subseteq T^*X$ and
$\Sigma \subseteq T^*Y$ closed conic subanalytic singular isotropics.
Then the K\"unneth formula holds:
\[
	\Sh_\Lambda(X) \otimes \Sh_\Sigma(Y) \simeq \Sh_{\Lambda \times \Sigma}(X \times Y).
\]
The dual of $\Sh_\Lambda(X)$ is $\Sh_{-\Lambda}(X)$.
And the equivalence
\[
	\Sh_{-\Lambda\times \Sigma}(X\times Y) \simeq \Sh_\Lambda(X)^\vee \otimes \Sh_\Sigma(Y) \simeq \Fun^L(\Sh_\Lambda(X), \Sh_\Sigma(Y))
\]
is given by the assignment
\[
	K \mapsto (-) \ast K.
\]
\end{recollection}
\begin{thm}\label{main-thm}
Let $K \in \Sh_{-\Lambda\times \Sigma}(X\times Y)$ be a sheaf kernel.
The convolution functor
\[
- \ast K\colon \Sh_\Lambda(X) \rightarrow \Sh_\Sigma(Y)
\]
preserves compact objects if and only if for every $x \in X$,
the restriction $K|_{\{x\} \times Y} \in \Sh_{\Sigma}(Y)$ is a compact object.
\end{thm}
\begin{proof}
By definition, the convolution preserves compact objects precisely
if $K \in \Sh_{-\Lambda \times \Sigma}(X \times Y) \simeq \Sh_{-\Lambda}(X) \otimes \Sh_\Sigma(Y)$ is a $\Sh_\Sigma(Y)$-proper object.
By \cref{v-proper-object-in-shv},
this is equivalent to
$K$ having compact stalks when viewed as a sheaf on $X$ valued in $\Sh_\Sigma(Y)$.
For each $x \in X$, the stalk functor and restriction to $\{x\} \times Y$ after
the K\"unneth equivalence both send exterior products $F \boxtimes G$ to
$F_x \otimes G$ in $\Sh(Y)$. Both functors preserve colimits, as does the
inclusion $\Sh_\Sigma(Y) \hookrightarrow \Sh(Y)$. The universal property of
the tensor product therefore identifies the stalk with $K|_{\{x\} \times Y}$
as an object of $\Sh_\Sigma(Y)$.
\end{proof}
\begin{cor}\label{perfect-stalk-and-compact-support}
Let $K \in \Sh_{-\Lambda \times \Sigma}(X \times Y)$ have perfect stalks.
If its restriction to $\{x\} \times Y$ has compact support for every $x \in X$,
then convolution with $K$ preserves compact objects.
\end{cor}
\begin{proof}
Combine \cref{compacts-in-shv} and \cref{main-thm}.
\end{proof}

\BiblatexSplitbibDefernumbersWarningOff
\printbibliography[keyword=alpha]
\printbibliography[notkeyword=alpha, heading=none]

\end{document}